\magnification\magstep1
\baselineskip14pt
\parskip3pt

\def\bib{\par\noindent\hangindent 20pt}
\def\bibbak{\kern-.5em}
\def\EGO{1}
\def\FGT{2}
\def\FKP{3}
\def\GS{4}
\def\GJ{5}
\def\Hi{6}
\def\Ki{7}
\def\Perm{8}
\def\Motz{9}
\def\PP{10}
\def\SR{11}
\def\MW{12}

\centerline{\bf Bracket notation for the `coefficient of' operator}
\smallskip
\centerline{Donald E. Knuth, Stanford University}
\bigskip\noindent
\centerline{[Corrected copy, 19 August 1993]}

\bigskip\noindent
When $G(z)$ is a power series in $z$, many authors now write
`$[z^n]\,G(z)$' for the coefficient of~$z^n$ in $G(z)$, using a
notation introduced by Goulden and Jackson in [\GJ, p.~1]. More
controversial, however, is the proposal of the same authors [\GJ,
p.~160] to let `$[z^n\!/n!]\,G(z)$' denote the coefficient of
$z^n\!/n!$, i.e., $n!$~times the coefficient of~$z^n$. An alternative
generalization of $[z^n]\,G(z)$, in which we define $[F(z)]\,G(z)$ to be a
linear function of both~$F$ and~$G$, seems to be more useful because
it facilitates algebraic manipulations. The purpose of this paper is
to explore some of the properties of such a definition. The remarks are
dedicated to Tony Hoare because of his lifelong interest in the improvement
of notations that facilitate manipulation.

\bigskip\noindent
{\bf Informal introduction.}\quad
In this paper `$[z^2+2z^3]\,G(z)$' will stand for the coefficient
of~$z^2$ plus twice the coefficient of~$z^3$ in~$G(z)$, when $G(z)$ is
a function of~$z$ for which such coefficients are well defined. More
generally, if $F(z)=f_0+f_1z+f_2z^2+\cdots$ and
$G(z)=g_0+g_1z+g_2z^2+\cdots\,$, we will let
$$[F(z)]\,G(z)=f_0g_0+f_1g_1+f_2g_2+\cdots$$
be the ``dot product'' of the vectors $(f_0,f_1,f_2,\ldots)$ and
$(g_0,g_1,g_2,\ldots)$, assuming that the infinite sum exists. Still
more generally, if $F(z)=\cdots
+f_{-2}z^{-2}+f_{-1}z^-+f_0+f_1z+f_2z^2+\cdots$ and $G(z)=\cdots
+g_{-2}z^{-2}+g_{-1}z^-+g_0+g_1z+g_2z^2+\cdots$ are doubly infinite
series, we will write
$$[F(z)]\,G(z)=\cdots
+f_{-2}g_{-2}+f_{-1}g_{-1}+f_0g_0+f_1g_1+f_2g_2+\cdots\;, \eqno(1)$$
again assuming convergence. (It is convenient to write `$z^-$' for
$1/z$, as in [\Perm].) 
The right side of~(1) is symmetric in~$F$ and~$G$, so we have a
commutative law:
$$[F(z)]\,G(z)=[G(z)]\,F(z)\,.\eqno(2)$$
There also is symmetry between positive and negative powers:
$$[F(z)]\,G(z)=[F(z^-)]\,G(z^-)\,.\eqno(3)$$

In particular, we will write $[1]\,G(z)$ for the constant term~$g_0$ of
a given doubly infinite power series $G(z)=\sum_ng_nz^n$. Notice that
$[z^n]\,G(z)=[1]\,z^{-n}G(z)$ and in fact
$$[F(z)]\,G(z)=[1]\,F(z^-)G(z)\,,\eqno(4)$$
when the product of series is defined in the usual way:
$$\sum_nh_nz^n=\biggl(\sum_nf_nz^n\biggr)\biggl(\sum_ng_nz^n\biggr)
\quad \Longleftrightarrow\quad h_n=\sum_{j+k=n}f_jg_k\,.\eqno(5)$$

Relation (4) gives us a useful rule for moving factors in and out of
brackets:
$$[F(z)]\,G(z)H(z)=[F(z)G(z^-)]\,H(z)\,.\eqno(6)$$
Both sides reduce to $[1]\,F(z^-)G(z)H(z)$, so they must be equal. This
rule is most often applied in a simple form such as
$$[z^n]\,z^3H(z)=[z^{n-3}]\,H(z)\,,$$
but it is helpful to remember the general principle (6). Similarly,
$$[F(z)G(z)]\,H(z)=[F(z)]\,G(z^-)H(z)\,.\eqno(7)$$

\bigskip\noindent
{\bf A paradox.}\quad
So far the extended bracket notation seems straightforward and
innocuous, but if we start to play with it in an undisciplined fashion
we can easily get into trouble. For example, one of the first uses we
might wish to make of relation~(1) is
$$\left[{z^n\over
1-z}\right]\,G(z)=g_n+g_{n+1}+g_{n+2}+\cdots\;,\eqno(8)$$ 
because $z^n\!/(1-z)=z^n+z^{n+1}+z^{n+2}+\cdots\,$. This,
unfortunately, turns out to be dangerous, if not outright fallacious. 

The danger is sometimes muted and we might be lucky. For example, if we
try combining (8) with (7) in the case $G(z)=1/(1-z)$ and
$H(z)=(1-z)^2=1-2z+z^2$, we get
$$\left[{z^n\over 1-z}\right]\,(1-z)^2=[z^n]\,{(1-z)^2\over
1-z^-}=[z^n]\,(z^2-z)\,.\eqno(9)$$ 
Sure enough, the sum $h_n+h_{n+1}+h_{n+2}+\cdots$ is nonzero in this
case only when $n=2$ and $n=1$, and (9) gives the correct answer. So
far so good.

But (7) and (8) lead to a contradiction when we apply them to the
trivial case $F(z)=H(z)=1$ and $G(z)=1/(1-z)$:
$$1=\left[{1\over 1-z}\right]\,1=[1]\,{1\over 1-z^-}=[1]\,{-z\over
1-z}=0\,.\eqno(10)$$
What went wrong?

\bigskip\noindent
{\bf Formal analysis.}\quad
To understand the root of the paradox (10), and to learn when (6) and
(7) are indeed valid rules of transformation, we need to know the
basic properties of double power series $\sum_ng_nz^n$. The general
theory can be found in Henrici [\Hi, \S4.4]; we will merely sketch it
here.

If $G(z)$ is analytic in an annulus $\alpha <\vert z\vert <\beta$, it
has a unique double series representation $G(z)=\sum_ng_nz^n$.
Conversely, any double power series that converges in an annulus
defines an analytic function there. The proof is based on the contour
integral formula 
$$G(z)={1\over 2\pi i}\;\oint_{\vert t\vert=\beta'}\;{G(t)\,dt\over
(t-z)} - {1\over 2\pi i}\;\oint_{\vert t\vert=\alpha'}\;{G(t)\,dt\over
(t-z)}\;, \eqno(11)$$
where $\alpha'$ is between $\vert z\vert$ and $\alpha$ while $\beta'$
is between~$\vert z\vert$ and~$\beta$. The quantity $1/(t-z)$ can be
expanded as $t^-(1+z/t+z^2\!/t^2+\cdots\,)$ when $\vert t\vert >\vert
z\vert$ and as $-z^-(1+t/z+t^2\!/z^2+\cdots\,)$ when $\vert t\vert
<\vert z\vert$. 

If $F(z)$ and $G(z)$ are both analytic for $\alpha <\vert z\vert
<\beta$, their product $H(z)$ is an analytic function whose
coefficients are given by~(5). Moreover, the infinite sum over all~$j$
and~$k$ with $j+k=n$ in~(5) is absolutely convergent: The terms are
$O\bigl((\alpha'/\beta')^k\bigr)$ as $k\rightarrow +\infty$ and
$O\bigl((\beta'/\alpha')^k\bigr)$ as $k\rightarrow -\infty$.

The coefficients of $G(z)$ in its double power series depend
on~$\alpha$ and~$\beta$. For example, suppose $G(z)=1/(2-z)$; we have
$${1\over 2-z}=\left\{\vcenter{\halign{$#$\hfil\quad&#\hfil\cr
{1\over 2}+{1\over 4}z+{1\over 8}z^2+\cdots\;,&when $\vert z\vert
<2$;\cr
\noalign{\medskip}
-z^--2z^{-2}-4z^{-3}-\cdots\;,&when $\vert z\vert >2$.\cr}}\right.
\eqno(12)$$
Thus if $F(z)=1/(2-z)+1/(2-z^-)$, there are three expansions
$$F(z)=\left\{\vcenter{\halign{$#$\hfil\quad&#\hfil\cr
{1\over 2}+\bigl({1\over 4}-1\bigr)z+\bigl({1\over 8}-2\bigr)
z^2+\bigl({1\over 16}-4\bigr)z^3+\cdots\;,&$\vert z\vert <{1\over 2}$;\cr
\noalign{\medskip}
\cdots +{1\over 8}z^{-2}+{1\over 4}z^-+1+{1\over 4}z+{1\over
8}z^2+\cdots \,,&${1\over 2}<\vert z\vert <2$;\cr
\noalign{\medskip}
\cdots +\bigl({1\over 16}-4\bigr)z^{-3}+\bigl({1\over 8}-2\bigr)z^{-2}
+\bigl({1\over 4}-1\bigr)z^-+{1\over 2}\,,&$\vert z\vert >2\,.$\cr}}\right.
\eqno(13)$$

Here's another example, this time involving a 
function that has an essential singularity instead of a pole:
$$e^{z/(1-z)}=\left\{\vcenter{\halign{$#$\hfil\quad&#\hfil\cr
1+z+{3\over 2}z^2+{13\over 6}z^3+{73\over 24}z^4+\cdots\,,
&$\vert z\vert <1\,;$\cr
\noalign{\medskip}
e^--e^-z^--{e^-\over 2}\,z^{-2}-{e^-\over 6}\,z^{-3} 
+{e^-\over 24}\,z^{-4}+\cdots\,,&$\vert z\vert >1\,.$\cr}}\right.
\eqno(14)$$
The coefficients when $\vert z\vert <1$ are $P_n/n!$, where $P_n$ is
the number of ``sets of lists'' of order~$n$ [\Motz].

\bigskip\noindent
{\bf Explaining the paradox.}\quad
The dependency of coefficients on $\alpha$ and $\beta$ makes our
notation $[F(z)]\,G(z)$ ambiguous; that is why we ran into trouble in the
paradoxical ``equation''~(10). We can legitimately use bracket
notation only when the context specifies a family  of
``safe'' functions---functions with well defined coefficients. 

The basic definition of $[F(z)]\,G(z)$ in (4) should be used only if the
product $F(z^-)G(z)$ is safe. Operation~(6), which moves a factor
$G(z)$  into the bracket, should be used only if
$F(z^-)G(z)H(z)$ is safe. Operation~(7), which removes a factor $G(z)$
from the bracket, should be used only if $F(z^-)G(z^-)H(z)$
is safe.

The root of our problem in (10) begins in (8), where we used the
expansion $F(z)=z^n\!/(1-z)=z^n+z^{n+1}+z^{n+2}+\cdots\;$; in other
words, $F(z^-)=z^{-n}\!/(1-z^-)=z^{-n}+z^{-n-1}+z^{-n-2}+\cdots\;$.
The latter expansion is valid only when $\vert z\vert >1$, so the
bracket notation of~(8) refers to coefficients in the region $1<\vert
z\vert <\infty$. In the last step of~(10), however, we said that
$[1]\,\bigl(-z/(1-z)\bigr)=0$, using coefficients from the region $\vert
z\vert <1$. The correct result for $\vert z\vert >1$~is
$$[1]\,{-z\over 1-z}=[1]\,{1\over 1-z^-}=[1]\,(\,\cdots
+z^{-2}+z^-+1)=1\,.$$

Bracket notation is most often used when $\vert z\vert$ is small, so we
should actually forget the ``rightward sum'' appearing in
equation~(8); it hardly ever yields the formula we want. The
``leftward sum'' rule 
$$\left[{z^n\over 1-z^-}\right]\,G(z)=\cdots
+g_{n-2}+g_{n-1}+g_n\eqno(15)$$ 
should be used instead, because
$z^n\!/(1-z^-)=\cdots +z^{n-2}+z^{n-1}+z^n$ is valid for $\vert
z^-\vert <1$. {\sl When the bracket notation $[F(z)]\,G(z)$ is being
used in the annulus $(\alpha,\beta)$, the functions $F(z^-)$ and
$G(z)$ should be analytic in $(\alpha,\beta)$.}
Note that $f(z^-)$ is analytic in $(\alpha,\beta)$ if and only if
$f(z)$ is analytic in $(\beta^-,\alpha^-)$.

\bigskip\noindent
{\bf Formal series.}\quad
Manipulations of generating functions are often done on formal power
series, when the coefficients are arbitrary and convergence is
disregarded. However, formal power series are not allowed to be
infinite in both directions; a~formal series $G(z)=\sum_ng_nz^n$ is
generally required to be a ``formal Laurent series''---a~series in
which $g_n=0$ for all sufficiently negative values of~$n$. We shall
call such series $L$-series for short. Similarly, we shall say that a
reverse formal Laurent series, in which $g_n=0$ for all sufficiently
{\it positive\/} values of~$n$, is an $R$-series. A~power series is
both an $L$-series and an $R$-series if and only if it is a polynomial
in~$z$ and~$z^-$.

Henrici [\Hi, \S1.2--1.8] shows that the normal operations on power
series---addition, subtraction, multiplication, division by nonzero,
differentiation, composition---can all be done rigorously on
$L$-series without regard to convergence. Thus $L$-series are ``safe''
functions: {\sl We can define bracket notation $[F(z)]\,G(z)$ by
rule~(4) whenever $F(z)$ is an $R$-series and $G(z)$ is an
$L$-series}.
Convergence is not then an issue. This definition provides the default
meaning of bracket notation, whenever no other context is specified.
The transformations in~(7) and~(8) are valid when the functions inside
brackets are $R$-series and the functions outside brackets are
$L$-series. Equations~(2) and~(3)
 should not be used unless $F$ and~$G$ are
both $L$-series and $R$-series.

In such cases paradoxes do not rear their ugly heads. The ill-fated
equation~(8) may fail, but equation~(15) is always true. 

\bigskip\noindent
{\bf Additional properties.}\quad
The bracket notation satisfies several identities in addition to~(2),
(3), (6), and~(7), hence we can often transform formulas in which it
appears. In the first place, the operation is linear in both operands:
$$\eqalignno{[a\,F(z)+b\,G(z)]\,H(z)&=a[F(z)]\,H(z)+b[G(z)]\,H(z)\,;&(16)\cr
\noalign{\smallskip}
[F(z)]\bigl(a\,G(z)+b\,H(z)\bigr)
&=a[F(z)]\,G(z)+b[F(z)]\,H(z)\,.&(17)\cr}$$

In the second place, there is a general multiplication law
$$[F_1(z)\,F_2(z)]\,G_1(z)\,G_2(z)=\sum_k\bigl([F_1(z)z^k]\,G_1(z)\bigr)
\bigl([F_2(z)z^{-k}]\,G_2(z)\bigr)\,.\eqno(18)$$
If $F_1(z)=F_2(z)=1$, this equation is simply the special case $n=0$
of~(5), and for general $F_1$ and~$F_2$ it follows from the special
case because we can replace $G_1(z)$ and $G_2(z)$ by $F_1(z^-)G_1(z)$
and $F_2(z^-)G_2(z)$ using~(7).

We also have
$$[F(z^m)]\,G(z^m)=[F(z)]\,G(z)\eqno(19)$$
for any nonzero integer~$m$; this equation, which includes (3) as the
special case $m=-1$, follows immediately from~(4) because
$[1]\,H(z)=[1]\,H(z^m)$. Equation~(19) suggests that we generalize bracket
notation to functions that are sums over nonintegral powers, in which
case $m$ would not need to be an integer. Then we could write (19) as
$$[F(z)]\,G(z^m)=[F(z^{1/m})]\,G(z)\,,\quad m\neq 0\,.\eqno(19')$$
Such generalizations, extending perhaps to integrals as well as to
sums, may prove to be quite interesting, but they will not be pursued
further here.

If $a$ is any nonzero constant, we have $[1]\,H(az)=[1]\,H(z)$. This rule
implies that
$[1]\,F(z^-)G(az)=[1]\,F(az^-)G(z)$, and (4) yields 
$$[F(z)]\,G(az)=[F(az)]\,G(z)\,.\eqno(20)$$
The special case where $F(z)$ is simply $z^m$ is, of course, already
familiar:
$$[z^m]\,G(az)=[(az)^m]\,G(z)=a^m[z^m]\,G(z)\,.$$

Bracket notation also interacts with differentiation in interesting
ways. We have, for instance,
$$[z^-]\,G'(z)=0\eqno(21)$$
for any function $G(z)=\sum_{n=-\infty}^{\infty}g_nz^n$. More
significantly,
$$[F(z)]\,z\,G'(z)=[z\,F'(z)]\,G(z)\,.\eqno(22)$$
Equation (21) is essentially the special case $F(z)=1$ of (22), but we
can also derive (22) from (21): Let 
$H(z)=F(z^-)G(z)$; then $0=[1]\,z\,H'(z)=[1]\,z\,\bigl(F(z^-)G'(z)
-z^{-2}F'(z)G(z)\bigr)$, hence
$[1]\,F(z^-)\,z\,G'(z)=[1]\,z^-F'(z^-)G(z)$, which is~(22).

Let $\vartheta$ be the operator $z\,{d\over dz}$. Then (22) implies by
induction on~$m$ that
$$[F(z)]\,\vartheta^m\,G(z)=[\vartheta^m\,F(z)]\,G(z)$$
for all integers $m\geq 0$, and we have
$$[F(z)]\,P(\vartheta)\,G(z)=[P(\vartheta)\,F(z)]\,G(z)\eqno(23)$$
for any polynomial $P$. If $F(z)=\sum_nf_nz^n$ and
$G(z)=\sum_ng_nz^n$, both sides of (23) evaluate to $\sum_n
P(n)f_ng_n$.

\bigskip\noindent
{\bf Additional variables.}\quad
When $G(w,z)$ is a bivariate generating function we also wish to write
$[w^mz^n]$ for the coefficient of $w^mz^n$ in~$G$. In general we can
define
$$[F(w,z)]\,G(w,z)=[1]\,F(w^-,z^-)\,G(w,z)\,,\eqno(24)$$
extending (4).

Variables must be clearly distinguished from constants. If $w$ and~$z$
are both variables, we have for instance $[z]\,wz=0$, while if $w$ is
constant we have $[z]\,wz=w$. If the set of variables is not clear from
the context, we can specify it by writing its elements as subscripts
on the brackets. For example,
$$[F(w)\,G(z)]_{w,z}\,H(w,z)=[G(z)]_z\bigl([F(w)]_w\,H(w,z)\bigr)\eqno(25)$$
because the former is $[1]_{w,z}\,F(w^-)G(z^-)H(w,z)$ while the latter
is
$$[1]_z\bigl(G(z^-)\,[1]_w\bigl(F(w^-)H(w,z)\bigr)\bigr)=[1]_z[1]_w
\,G(z^-)\,F(w^-)\,H(w,z)$$ 
and $[1]_{w,z}=[1]_w[1]_z$. 

After we have evaluated the parenthesis on the right side of~(25), the
ambiguity disappears, because $w$ is no longer present. For example,
if $m\geq 0$ we have
$$[w^mz^n]\,{1\over 1-wF(z)}=[z^n]\left([w^m]_w\,{1\over 1-wF(z)}\right)
=[z^n]\,F(z)^m\,,\eqno(26)$$
where brackets without subscripts assume that both $w$ and $z$ are
variables. Similarly
$$\eqalignno{[w^mz^n]\,e^{wF(z)}&=[z^n]\,{F(z)^m\over m!}\,,&(27)\cr
\noalign{\smallskip}
[w^mz^n]\,G\bigl(w\,F(z)\bigr)\,H(z)&=[z^n]\,F(z)^m\,H(z)\,[w^m]\,
G(w)\,.&(28)\cr}$$

Suppose $w$ and $z$ are variables. Then laws (19) and (20) extend to
$$\eqalignno{[F(w,z)]\,G(aw,z)&=[F(aw,z)]\,G(w,z)\,,\ \qquad a\neq
0\,;&(29)\cr
\noalign{\smallskip}
[F(w^m,z)]\,G(w^m,z)&=[F(w,z)]\,G(w,z)\,,\quad\qquad{\rm integer}\ m\neq
0\,;&(30)\cr
\noalign{\smallskip}
[F(w,w^mz)]\,G(w,w^mz)&=[F(w,z)]\,G(w,z)\,;&(31)\cr}$$
and we have indeed the general rule
$$[F(a^-w^kz^l,b^-w^mz^n)]\,G(aw^kz^l,bw^mz^n)=[F(w,z)]\,G(w,z)\eqno(32)$$
when $a\neq 0$, $b\neq 0$, and 
$\left\vert{k\ l\atop m\ n}\right\vert\neq 0$, i.e., 
$kn\neq lm$. A~similar formula applies
with respect to any number of variables.

The following example from the theory of random graphs [\FKP, (10.10)
and (10.14)] illustrates how these rules are typically applied.
Suppose we want to evaluate the coefficient of $[w^mz^n]$ in the
expression $e^{U(wz)/w+V(wz)}$, where $U$ and~$V$ are known functions
with $U(0)=0$.
The two-variable problem is reduced to a one-variable problem as
follows:
$$\eqalignno{[w^mz^n]\,e^{U(wz)/w+V(wz)}
&=[(w^-)^{n-m}(wz)^n]\,e^{U(wz)w^-+V(wz)}\cr
\noalign{\smallskip}
&=[w^{n-m}z^n]\,e^{U(z)w+V(z)}\cr
\noalign{\smallskip}
&={1\over (n-m)!}\,[z^n]\,U(z)^{n-m}e^{V(z)}\,,&(33)\cr}$$
by (32) with $F(w,z)=w^{n-m}z^n$, $G(w,z)=e^{U(z)w+V(z)}$, $a=b=1$,
$k=-1$, $l=0$, $m=n=1$. The final step uses (28) with $F(z)=U(z)$,
$G(w)=e^w$, and $H(z)=e^{V(z)}$.

As before, we need to check that the functions are safe before we can
guarantee that such manipulations are legitimate. For formal power
series, the functions inside brackets should be $R$-series and the
functions outside should be $L$-series. This condition holds in each
step of~(33) because $U(0)=0$.

\bigskip\noindent
{\bf Additional identities.}\quad
The bracket notation also obeys more complex laws that deserve further
study. For example, Gessel and Stanton [\GS, Eq.~(3)] have shown among
other things that
$$[F(w,z)]\,{G(w,z)\over
1-wz}=\bigl[F\bigl(w(1+z^-),z(1+w^-)\bigr)\bigr]\,G\left({w\over
1+z}\,,\,{z\over 1+w}\right)\,.\eqno(34)$$
If we set $F(w,z)=w^kz^l$ and $G(w,z)=(1+w)^m(1+z)^n/(1-wz)^{m+n}$,
Gessel and Stanton observe that  we
obtain Saalsch\"utz's identity after some remarkable cancellation:
$$\eqalignno{\sum_r\,{m\choose k-r}{n\choose l-r}{m+n+r\choose r}
&=[w^k(1+z^-)^kz^l(1+w^-)^l]\,(1+w)^m(1+z)^n\cr
\noalign{\medskip}
&=[w^kz^l]\,(1+w)^{m+l}(1+z)^{n+k}\cr
\noalign{\medskip}
&={m+l\choose k}{n+k\choose l}\,.&(35)\cr}$$ 
And if we set $F(w,z)=w^{l+n}z^{m+n}$,
$G(w,z)=(w-z)^{l+m}\!/(1-wz)^{l+m}$, the left side of (34) reduces to
$$[w^{l+n}z^{m+n}]\,{(w-z)^{l+m}\over (1-wz)^{l+m+1}}=(-1)^m\,
{(l+m+n)!\over l!\,m!\,n!}\;;\eqno(36)$$
the right side is
$$\eqalignno{&[w^{l+n}(1+z^-)^{l+n}z^{m+n}(1+w^-)^{m+n}]\,(w-z)^{l+m}\cr
\noalign{\smallskip}
&\qquad\qquad =[w^{l+n}z^{m+n}]\,(w-z)^{l+m}(1+w)^{m+n}(1+z)^{l+n}\cr
\noalign{\smallskip}
&\qquad\qquad =\sum_k(-1)^{k+m}{l+m\choose k+m}{m+n\choose
k+n}{n+l\choose k+l}\,.&(37)\cr}$$
The fact that $(36) =(37)$ is Dixon's identity [\Ki, exercise
1.2.6--62].

Equation (34) can be generalized to $n$ variables, and we can replace
the~`1' on the right by any nonzero constant~$a$:
$$\eqalignno{[F(z_1,&\ldots,z_n)]\,{G(z_1,\ldots,z_n)\over
1-z_1\,\ldots\,z_n}\cr
\noalign{\medskip}
&=\bigl[F\bigl(z_1(a+z_2^-),\ldots,z_n(a+z_1^-)\bigr)\bigr]\,
G\left({z_1\over a+z_2}\,,\ldots,\,{z_n\over
a+z_1}\right)\,.&(38)\cr}$$
It suffices to prove this when $F(z_1,\ldots,z_n)=1$ and
$G(z_1,\ldots,z_n)=z_1^{m_1}\,\ldots\,z_n^{m_n}$, in which case both
sides are~0 unless $m_1=\cdots =m_n\leq 0$, when both sides are~1.
Equation (38) holds in particular when $n=1$:
$$[F(z)]\,{G(z)\over 1-z}=[F(1+az)]\,G\left({z\over
a+z}\right)\,,\qquad a\neq 0\,.\eqno(39)$$

Returning to the case of a single variable, we should also state the
general rule for composition of series:
$$G\bigl(F(z)\bigr)=\sum_nF(z)^n\,[z^n]\,G(z)\,.\eqno(40)$$
Special conditions are needed to ensure that this infinite sum is well
defined.

\bigskip\noindent
{\bf Lagrange's inversion formula.}\quad
Let $F(z)=f_1z+f_2z^2+f_4z^3+\cdots\;$, with $f_1\neq 0$, and let
$G(z)$ be the inverse function so that
$$F\bigl(G(z)\bigr)=G\bigl(F(z)\bigr)=z\,.\eqno(41)$$
Lagrange's celebrated formula for the coefficients of $G$ can be
expressed in bracket notation in several ways; for example, we have
$$n[z^n]\,G(z)^m=m\,[z^{-m}]\,F(z)^{-n}\,,\eqno(42)$$
for all integers $m$ and $n$.

One way to derive (42), following Paule [\PP], is to note first that
(40) implies
$$z^m=G\bigl(F(z)\bigr)^m=\sum_k F(z)^k\,[z^k]\,G(z)^m.\eqno(43)$$
Differentiating with the $\vartheta$ operator and dividing by $F(z)^n$ yields
$${mz^m\over F(z)^n}=\sum_k kF(z)^{k-1-n}\vartheta F(z)\,[z^k]\,G(z)^m.
\eqno(44)$$
Now we will study the constant terms of (44). If $k\ne n$,
$$[1]\,F(z)^{k-1-n}\vartheta F(z) = [1]\,{\vartheta\bigl(F(z)^{k-n}\bigr)\over
k-n}=0\,,\eqno(45)$$
by (22). And if $k=n$,
$$[1]\,{\vartheta F(z)\over F(z)}=
[1]\,{f_1+2f_2z+3f_3z^2+\cdots\over f_1+f_2z+f_3z^2+\cdots}=1\,,\eqno(46)$$
because $f_1\ne0$. Therefore the constant terms of (44) are
$$[1]\,{mz^m\over F(z)^n}=n\,[z^n]\,G(z)^m\,;$$
this is Lagrange's formula (42).

\bigskip\noindent
{\bf Conclusions.}\quad
Many years of experience have confirmed the great importance of
generating functions in the analysis of algorithms, and we can
reasonably expect that some fluency in manipulating the
``coefficient-of'' operator will therefore be rewarding.

If, for example, we are faced with the task of simplifying a formula
such as
$$\sum_k{m\choose k}[z^{n-k}]\,F(z)^k\,,$$
a rudimentary acquaintance with the properties of brackets will tell
us that it can be written as $\sum_k{m\choose k}[z^n]\,z^k\,F(z)^k$
and then summed to yield
$$[z^n]\bigl(1+z\,F(z)\bigr)^m\,.$$
We have seen several examples above in which formulas that are far
less obvious can be derived rapidly by bracket manipulation, when we
use quantities more general than monomials inside the brackets.

In most applications we use bracket notation in connection with formal
Laurent series, in which case it is important to remember that our
identities for $[F(z)]\,G(z)$ require $G(z)$ to have only finitely
many {\it negative\/} powers of~$z$ while $F(z)$ must have only
finitely many {\it positive\/} powers. If we write, for example,
$$\left[{z^n\over z-1}\right]\,G(z)\,,\eqno(47)$$
we should think of the quantity in brackets as an infinite series
$$z^{n-1}+z^{n-2}+z^{n-3}+\cdots$$
that descends to arbitrarily {\it negative\/} powers of $z$; the
bracket notation then denotes the sum
$g_{n-1}+g_{n-2}+g_{n-3}+\cdots\;$, which will be finite. We have seen
that other interpretations of bracket notation are possible for
functions analytic in an annulus; but great care must be taken to
avoid paradoxes in such cases, hence the extra effort might not be
worthwhile.

Bracket notation, like all notations, is ``dispensable,'' in the sense
that we can prove the same theorems without it as with it. But the use
of a good notation can shorten proofs and help us see patterns that
would otherwise be difficult to perceive.

Let us close with one more example, illustrating that the notation
(47) helps to simplify some of the formulas in [\FGT]. The {\it coupon
collector's problem\/} asks for the expected number of trials needed
to obtain $n$~distinct coupons from a set~$C$ of~$m$ given coupons,
where each trial independently produces coupon~$c$ with probability
$p(c)$. Theorem~2 of [\FGT] says, when rewritten in the notation
discussed above, that this expected number is
$$\int_0^{\infty}\,\left[{z^n\over z-1}\right]\,\prod_{c\in C}\,
\bigl(1+z(e^{p(c)t}-1)\bigr)\,e^{-t}\,dt\,.\eqno(48)$$
We can evaluate (48) by expanding the integrand as follows:
$$\eqalign{\left[{z^n\over z-1}\right]\,\prod_{c\in C}\,
\bigl(1+z(e^{p(c)t}-1)\bigr)
&=\sum_{\scriptstyle {B\subseteq C}\atop \scriptstyle {\vert B\vert <n}}\;
\prod_{c\in B}\,(e^{p(c)t}-1)\cr
\noalign{\medskip}
&=\sum_{\scriptstyle {A\subseteq B\subseteq C}\atop\scriptstyle {\vert
B\vert <n}}\,(-1)^{\vert B\vert-\vert A\vert}e^{p(A)t}\cr
\noalign{\medskip}
&=\sum_{\scriptstyle {A\subseteq C}\atop\scriptstyle {\vert A\vert <n}}
e^{p(A)t}\sum_{\vert A\vert \leq k<n}(-1)^{k-\vert A\vert}
{\vert C\vert-\vert A\vert\choose k-\vert A\vert}\cr
\noalign{\medskip}
&=\sum_{\scriptstyle{A\subseteq C}\atop\scriptstyle{\vert A\vert <n}}
e^{p(A)t}(-1)^{n-1-\vert A\vert}
{\vert C\vert-\vert A\vert-1\choose\vert C\vert -n}\,,\cr}$$
where $p(A)$ denotes $\sum_{a\in A}p(a)$. The integral (48) therefore
is
$$\sum_{\scriptstyle{A\subseteq C}\atop\scriptstyle{\vert
A\vert <n}}(-1)^{n-1-\vert A\vert}\left.{\vert C\vert-\vert
A\vert -1\choose \vert C\vert
-n}\right/\bigl(1-p(A)\bigr)\,.\eqno(49)$$
(This is Corollary 3 of [\FGT], which was stated without proof.)

\bigskip\noindent
{\bf Related work.}\quad
Steven Roman's book on umbral calculus [\SR] develops extensive
properties of his notation $\langle G(t)\mid F(x)\rangle$, which
equals $\sum_{n\geq 0}f_ng_n$ when $F(x)=\sum_{n\geq 0}f_nx^n$ and
$G(t)=\sum_{n\geq 0}g_nt^n\!/n!$; the function $F(x)$ in these
formulas must be a polynomial. Thus, if $D$ is the operator $d/dx$,
Roman's $\langle G(t)\mid F(x)\rangle$ is the constant term of the
polynomial $G(D)\,F(x)$. Chapter~6 of [\SR] considers generalizations in
which $\langle G(t)\mid F(x)\rangle$ is defined to be $\sum_{n\geq
0}f_ng_n$ when $G(t)=\sum_{n\geq 0}g_nt^n\!/c_n$ and $c_n$ is an
arbitrary sequence of constants; the case $c_n=1$ corresponds to the
special case of bracket notation $[F(z)]\,G(z)$ when $F$ and~$G$
involve no negative powers of~$z$. Roman traces the theory back to a
paper by Morgan Ward [\MW].

G. P. Egorychev's book [\EGO] includes a great many examples that demonstrate
the value of coefficient extraction in the midst of formulas.

\bigskip\noindent
{\bf Open problems.}\quad
One reason formal power series are usually restricted to $L$-series is
that certain doubly infinite power series are divisions of zero. For
example, $\sum_{n=-\infty}^{\infty}z^n$ is a divisor of zero because
multiplication by $1-z$ annihilates it. (This series causes no problem
in the theory of non-formal power series because it does not converge
for any value of~$z$.) All double series having the form
$\sum_nn^m\alpha^nz^n$ for $\alpha\neq 0$ and integer $m\geq 0$ can
also be shown to be divisors of zero.   Question: Do there exist
divisors of zero besides finite linear combinations of the double
series just mentioned? Conjecture: There is no nonzero double series
$F(z)$ such that $e^zF(z)=0$. (A~counterexample would necessarily be
divergent.) 

It may be possible and interesting to extend the theory of formal
Laurent series to arbitrary functions of the form $F(z)\sum_ng_nz^n$,
where $g_n$ is zero for all sufficiently negative~$n$ and where $F(z)$
is analytic for $0<\vert z\vert <\infty$. 

\bigskip\noindent
{\bf Acknowledgments.}\quad
I wish to thank Edsger and Ria Dijkstra for the splendid opportunity to
write this paper in the guest room of their Texas home, and Peter Paule
for his penetrating comments on the first draft.

\vfill\eject

\bigskip
\centerline{Bibliography}

\medskip\bib
[\EGO]\quad
G. P. Egorychev, {\sl Integral Representation and the Computation of
Combinatorial Sums\/} (Providence, Rhode Island: American Mathematical
Society, 1984).

\medskip\bib
[\FGT]\quad
Philippe Flajolet, Dani\`ele Gardy, and Lo\"ys Thimonier, ``Birthday
paradox, coupon collectors, caching algorithms and self-organizing
search,'' {\sl Discrete Applied Mathematics\/ \bf 39} (1992), 207--229.

\medskip\bib
[\FKP]\quad
Philippe Flajolet, Donald E. Knuth, and Boris Pittel, ``The first
cycles in an evolving graph,'' {\sl Discrete Mathematics\/ \bf 75}
(1989), 167--215.

\medskip\bib
[\GS]\quad
Ira Gessel and Dennis Stanton, ``Short proofs of Saalsch\"utz's and
Dixon's theorems,'' {\sl Journal of Combinatorial Theory\/ \bf A38}
(1985), 87--90.

\medskip\bib
[\GJ]\quad
I. P. Goulden and D. M. Jackson, {\sl Combinatorial Enumeration\/}
(New York: Wiley, 1983).


\medskip\bib
[\Hi]\quad
Peter Henrici, {\sl Applied and Computational Complex Analysis},
Volume~1 (New York: Wiley, 1974).

\medskip\bib
[\Ki]\quad
Donald E. Knuth, {\sl Fundamental Algorithms}, Volume~1 of {\sl The
Art of Computer Programming\/} (Reading, Massachusetts:
Addison\kern.05em--Wesley, 1968).

\medskip\bib
[\Perm]\quad
Donald E. Knuth, ``Efficient representation of perm groups,''
{\sl Combinatorica\/ \bf 11}\allowbreak
 (1991), 33--43.

\medskip\bib
[\Motz]\quad
T. S. Motzkin, ``Sorting numbers for cylinders and other
classification numbers,'' {\sl Proceedings of Symposia in Pure
Mathematics\/ \bf 19} (1971), 167--176.

\medskip\bib\bibbak
[\PP]\quad
Peter Paule, ``Ein neuer Weg zur $q$-Lagrange Inversion,''
{\sl Bayreuther Mathematische Schriften\/ \bf18} (1985), 1--37.

\medskip\bib\bibbak
[\SR]\quad
Steven Roman, {\sl The Umbral Calculus\/} (Orlando, Florida: Academic
Press, 1984).

\medskip\bib\bibbak
[\MW]\quad
Morgan Ward, ``A calculus of sequences,'' {\sl American Journal of
Mathematics\/ \bf 58} (1936), 255--266.

\bye